# Some Characterizations of Weakly Pseudo Primary 2-Absorbing Submodules in Terms of some Types of Modules

Omar Hisham Taha[1]*, Marwa Abdullah Salih [2]

[1]Department of Mathematics, College of Computer Science and Mathematics, Tikrit University, Tikrit, Iraq

[2] Department of Mathematics, College of Education for Pure Sciences, Tikrit University, Tikrit, Iraq

**Abstract**

All rings are commutative with $1 \neq 0$, and all modules are unital. The purpose of this paper is to investigate the characterizations of weakly pseudo primary 2-absorbing sub-module in terms of some types of modules. We provide characterizations for the class of multiplication modules with the help of some types of modules such as faithful, non-singular, Z-regular, and projective modules. Furthermore, we add some conditions to prove the residual of a weakly pseudo primary 2-absorbing sub-module is a weakly pseudo primary 2-absorbing ideal.

## 1. Introduction

In this paper, the rings $R$ is a commutative ring with a non-zero identity, and $D$ is a unitary left $R$-module. During the last 13 years, the notion of 2-absorbing sub-modules and weakly 2-absorbing sub-modules has been previously investigated by [1], and its generalizations to 2-absorbing ideal, where 2-absorbing ideal was first introduced in 2007 by Badawi [2], respectively. A sub-module $A \subsetneq D$ of an $R$–module $D$ is called 2-absorbing (weakly 2-absorbing) if whenever $r_1 r_2 d \in A$ ($0 \neq r_1 r_2 d \in A$) for some $r_1, r_2 \in R$, $d \in D$, then $r_1 d \in A$ or $r_2 d \in A$ or $r_1 r_2 \in [A:_R D]$. Following that, Mostafanasab, Yetkin, Tekir and Daran [3] introduced the concept of primary- 2-absorbing sub-module and weakly primary 2-absorbing sub-module as generalizations of primary sub-module and weakly primary sub-module, respectively, a sub-module $A \subsetneq D$ of an $R$–module $D$ is called a primary-2-absorbing, if $r_1 r_2 d \in A$, for $r_1, r_2 \in R$, $d \in D$, implies that either $r_1 d \in rad_D(A)$ or $r_2 d \in rad_D(A)$ or $r_1 r_2 \in [A + Soc(D):_R D]$. Also, in [4] introduced the concept of pseudo primary 2-absorbing sub-module. A sub-module $A \subsetneq D$ of an $R$–module $D$, is called a pseudo primary 2-absorbing sub-module $D$, if $r_1 r_2 d \in A$, for $r_1, r_2 \in R$, $d \in D$, implies that either $r_1 d \in rad_D(A) + Soc(D)$ or $r_2 d \in rad_D(A) + Soc(D)$ or $r_1 r_2 \in [A + Soc(D):_R D]$. The concept of weakly pseudo primary 2-absorbing is a generalization of primary 2-absorbing and pseudo primary 2-absorbing sub-module in [5]. A sub-module $A \subsetneq D$ of an $R$–module $D$, is called a weakly pseudo primary 2-absorbing sub-module $D$, if $0 \neq r_1 r_2 d \in A$, for $r_1, r_2 \in R$, $d \in D$, implies that either $r_1 d \in rad_D(A) + Soc(D)$ or $r_2 d \in rad_D(A) + Soc(D)$ or $r_1 r_2 \in [A + Soc(D):_R D]$. An ideal $Q$ of a ring $R$ is said to be a weakly pseudo primary 2-absorbing ideal of $R$ if $Q$ is a weakly pseudo primary -2-absorbing $R$–submodule an $R$–module of $R$. An $R$-module $D$ is called multiplication if every sub-module $A$ of $D$ is of the from $A = QD$ for some ideal $Q$ of $R$. Equivalent to $A = [A:_R D]D$. It's well known that a cyclic module is multiplication [6]. An $R$-module $D$ is called said to be a Z-regular if for any $d \in D$ there exists $f \in D^* = Hom_R(D, R)$ such that $d = f(d)d$ [7]. Recall that an $R$-module $D$ is a projective if for all $R$-epimorphism $f: M_1 \to M_2$ and for every $R$-homomorphism $g: D \to M_2$, there exists an $R$-homomorphism $h: D \to M_1$ such that $f \circ h = g$ [8]. An $R$-module $D$ said to be non-singular if $Z(D) = 0$, where $Z(D) = \{d \in D: dQ = (0)\}$ for some ideal $Q$ of $R$ [9]. Also, an $R$-module $D$ is called faithful if $Ann_R(D) = (0)$ where $Ann(D) = \{r \in R: rD = (0)\}$. Recall that an $R$-module $D$ be a finitely generated if $D = (x_1, x_2 \ldots, x_n) = Rx_1 + Rx_2 + \ldots + Rx_n$, where $x_1, x_2, \ldots, x_n \in D$ [10]. In section two of this paper, we provide characterizations for many modules, including non-singular modules, multiplication modules, faithful finitely produced modules, projective modules, and Z-regular modules. And show under certain conditions, the residual of weakly pseudo primary 2-absorbing sub-modules is a weakly pseudo primary 2-absorbing ideal, see Proposition 2.12, Proposition 2.16, Proposition 2.19, and Proposition 2.22. Also, we show that under certain conditions, if $Q$ is weakly pseudo primary 2-absorbing ideal of $R$, then $QD$ is a weakly pseudo primary 2-absorbing sub-module of an $R$-module $D$ see Proposition 2.28, Proposition 2.30, and Proposition 2.32.

---

* omar.h.tahamm2314@st.tu.edu.iq

## Some Characterizations of Weakly Pseudo Primary 2-Absorbing Sub-module in Terms some Types of Modules

**Proposition 2.1**. [5] Let $A \subsetneq D$ be sub-module of $R$-module $D$. Then $A$ is a weakly pseudo primary 2-absorbing sub-module of an R-module D if and only if $(0) \neq r_1 r_2 T \subseteq A$ for $r_1 r_2 \in R$ and $T$ is sub-module of $D$, implies either $r_1 T \subseteq rad_D(A) + Soc(D)$ or $r_2 T \subseteq rad_D(A) + Soc(D)$ or $r_1 r_2 \in [A + Soc(D):_R D]$.

**Proposition 2.2**. Let $D$ be a multiplication $R$-module and $A \subsetneq D$ be a sub-module of $D$. Then $A$ is a weakly pseudo primary 2-absorbing sub-module of $D$ if and only if whenever $(0) \neq T_1 T_2 T_3 \subseteq A$, for some sub-modules $T_1, T_2, T_3$ of $D$, implies that either $T_1 T_3 \subseteq rad_D(A) + Soc(D)$ or $T_2 T_3 \subseteq rad_D(A) + Soc(D)$ or $T_1 T_2 \subseteq A + Soc(D)$.

**Proof**. Assume that $(0) \neq T_1 T_2 T_3 \subseteq A$, for some sub-modules $T_1, T_2, T_3$ of $D$. Since $D$ is a multiplication then $T_1 = Q_1 D, T_2 = Q_2 D, T_3 = Q_3 D$ for some ideals $Q_1, Q_2, Q_3$ in $R$, it follows that $(0) \neq Q_1 D Q_2 D Q_3 D = Q_1 Q_2(Q_3 D) \subseteq A$, since $A$ is a weakly pseudo primary 2-absorbing sub-module by Proposition 2.1 we obtain $Q_1 Q_3 D \subseteq rad_D(A) + Soc(D)$ or $Q_2 Q_3 D \subseteq rad_D(A) + Soc(D)$ or $Q_1 Q_2 \subseteq [A + Soc(D):_R D]$, which implies that either $(Q_1 D)(Q_3 D) \subseteq rad_D(A) + Soc(D)$ or $(Q_2 D)(Q_3 D) \subseteq rad_D(A) + Soc(D)$ or $(Q_1 D)(Q_2 D) \subseteq A + Soc(D)$, that is either $T_1 T_3 \subseteq rad_D(A) + Soc(D)$ or $T_2 T_3 \subseteq rad_D(A) + Soc(D)$ or $T_1 T_2 \subseteq A + Soc(D)$.

Conversely, $(0) \neq r_1 r_2 T \subseteq A$, for some sub-module $T$ of $D$ and $r_1, r_2 \in R$ since $D$ is a multiplication, then $T = QD$ for some ideal $Q$ in $R$, then $(0) \neq r_1 r_2 Q D \subseteq A$, it follows by hypothesis either $r_1 Q D \subseteq rad_D(A) + Soc(D)$ or $r_2 Q D \subseteq rad_D(A) + Soc(D)$ or $r_1 r_2 \in [A + Soc(D):_R D]$, that is $r_1 T \subseteq rad_D(A) + Soc(D)$ or $r_2 T \subseteq rad_D(A) + Soc(D)$ or $r_1 r_2 D \in A + Soc(D)$ by Proposition 2.1 implies that $A$ is a weakly pseudo primary 2-absorbing sub-module of $D$. □

The following corollary is obtained as an immediate outcome of Proposition 2.2.

**Corollary 2.3**. Let $D$ be a multiplication $R$-module and $A \subsetneq D$ be a sub-module of $D$. Then $A$ is a weakly pseudo primary 2-absorbing sub-module of $D$ if and only if whenever $(0) \neq T_1 T_2(d) \subseteq A$, for some sub-modules $T_1, T_2$ of $D$ and $d \in D$, implies that either $T_1(d) \subseteq rad_D(A) + Soc(D)$ or $T_2(d) \subseteq rad_D(A) + Soc(D)$ or $T_1 T_2 \subseteq A + Soc(D)$.

**Corollary 2.4**. Let $D$ be a cyclic R-module and $A \subsetneq D$ be a sub-module of $D$. Then $A$ is a weakly pseudo primary 2-absorbing sub-module of $D$ if and only if whenever $(0) \neq T_1 T_2 T_3 \subseteq A$, for some sub-modules $T_1, T_2 T_3$ of $D$, implies that either $T_1 T_3 \subseteq rad_D(A) + Soc(D)$ or $T_2 T_3 \subseteq rad_D(A) + Soc(D)$ or $T_1 T_2 \subseteq A + Soc(D)$.

**Corollary 2.5**. Let $D$ be a cyclic R-module and $A \subsetneq D$ be a sub-module of $D$. Then $A$ is a weakly pseudo primary 2-absorbing sub-module of $D$ if and only if whenever $(0) \neq T_1 T_2 d \subseteq A$, for some sub-modules $T_1, T_2$ of $D$ and $d \in D$, implies that either $T_1 d \subseteq rad_D(A) + Soc(D)$ or $T_2 d \subseteq rad_D(A) + Soc(D)$ or $T_1 T_2 \subseteq A + Soc(D)$.

**Corollary 2.6**. Let $D$ be a multiplication $R$-module and $A \subsetneq D$ be a sub-module of $D$. Then $A$ is a weakly pseudo primary 2-absorbing sub-module of $D$ if and only if whenever $(0) \neq T(d)D \subseteq A$, for some sub-module $T$ of $D$ and $d \in D$, implies that either $TD \subseteq rad_D(A) + Soc(D)$ or $(d)D \subseteq rad_D(A) + Soc(D)$ or $T(d) \subseteq A + Soc(D)$.

**Remark 2.7**. [5] The residual of a weakly pseudo primary 2-absorbing sub-module of an $R$-module $D$ is not necessarily to be a weakly pseudo primary 2-absorbing ideal of $R$. and this can be shown by the following example.

**Example 2.8**. Consider $A = (\overline{0})$ be a sub-module of $\mathbb{Z}$ - module $\mathbb{Z}_{30}$. As $A$ is a weakly pseudo primary 2-absorbing by detention, but $[A :_\mathbb{Z} \mathbb{Z}_{30}] = 30\mathbb{Z}$ is not weakly pseudo primary 2-absorbing ideal of $\mathbb{Z}$. Since $2.3.5 = 30 \in 30\mathbb{Z}$ for $2,3,5 \in \mathbb{Z}$. But $2.5 = 10 \notin rad_D(30\mathbb{Z}) + Soc(\mathbb{Z}) = 30\mathbb{Z} + (0) = 30\mathbb{Z}$ and $3.5 = 15 \notin rad_D(30\mathbb{Z}) + Soc(\mathbb{Z}) = 30\mathbb{Z}$ and $2.3 = 6 \notin [30\mathbb{Z} + Soc(\mathbb{Z}):_\mathbb{Z} \mathbb{Z}] = 30\mathbb{Z}$.

**Lemma 2.9**. [7] *If $D$ is Z-regular of an R-module then $Soc(D) = Soc(R)D$.*

**Lemma 2.10**. [11] Let $A \subsetneq D$ be a sub-module of multiplication $R$-module $D$. Then $rad_D(A) = \sqrt{[A:_R D]}D$.

**Proposition 2.11**. [5] Let $A$ be a weakly pseudo primary 2-absorbing sub-module of an $R$-module $D$ if and only if $(0) \neq Q_1 Q_2 T \subseteq A$ for some ideals $Q_1, Q_2$ *of $R$* and sub-module $T$ *of $D$*, implies that either $Q_1 T \subseteq rad_D(A) + Soc(D)$ or $Q_2 T \subseteq rad_D(A) + Soc(D)$ or $Q_1 Q_2 \subseteq [A + Soc(D):_R D]$.

Under certain conditions, the residual of a weakly pseudo primary 2-absorbing sub-module is a weakly pseudo primary 2-absorbing ideal, as shown in the following Propositions.

**Proposition 2.12**. Let $D$ be a Z-regular multiplication $R$-module and $A \subsetneq D$ be a sub-module of $D$. Then A is a weakly pseudo primary 2-absorbing sub-module of $D$ if and only if $[A:_R D]$ is a weakly pseudo primary 2-absorbing ideal of $R$.

**Proof**. Assume that $(0) \neq T_1 T_2 T_3 \subseteq [A:_R D]$ for some sub-module $T_1, T_2, T_3$ of $R$, it follows that $(0) \neq T_1 T_2 T_3 D \subseteq A$. But $D$ is a multiplication, then there exist three ideals $Q_1, Q_2, Q_3$ in $R$ such that $T_1 = Q_1D$, $T_2 = Q_2D$ and $T_3 = Q_3D$. Implies that $(0) \neq Q_1Q_2Q_3D \subseteq A$. Then by Proposition 2.2, either $Q_1Q_3D \subseteq rad_D(A) + Soc(D)$ or $Q_2Q_3D \subseteq rad_D(A) + Soc(D)$ or $Q_1Q_2D \subseteq A + Soc(D)$. Since $D$ is a Z-regular, then by Lemma 2.9, $(Soc(D) = Soc(R)D)$, and since $D$ is multiplication, then $A = [A:_R D]D$, and by Lemma 2.10 $rad_D(A) = \sqrt{[A:_R D]}D$, thus either $Q_1Q_3D \subseteq \sqrt{[A:_R D]}D + Soc(R)D$ or $Q_2Q_3D \subseteq \sqrt{[A:_R D]}D + Soc(R)D$ or $Q_1Q_2D \subseteq [A:_R D]D + Soc(R)D$. That is either $Q_1Q_3 \subseteq \sqrt{[A:_R D]} + Soc(R)$ or $Q_2Q_3 \subseteq \sqrt{[A:_R D]} + Soc(R)$ or $Q_1Q_2 \subseteq [A:_R D] + Soc(R) = [[A:_R D] + Sco(R):_R R]$. Therefore, by Proposition 2.11. $[A:_R D]$ is a weakly pseudo primary 2-absorbing ideal of $R$.

Conversely, Assume that $(0) \neq r_1r_2T \subseteq A$ for some sub-module $T$ of $R$, and $r_1, r_2 \in R$. But $D$ is a multiplication, then there exists an ideal $Q$ in $R$ such that $T = QD$. Implies that $(0) \neq r_1r_2QD \subseteq A$. That is $(0) \neq r_1r_2Q \in [A:_R D]$. Since $[A:_R D]$ is a weakly pseudo primary 2-absorbing ideal of $R$, Then by Proposition 2.1, either $r_1Q \subseteq \sqrt{[A:_R D]} + Soc(R)$ or $r_2Q \subseteq \sqrt{[A:_R D]} + Soc(R)$ or $r_1r_2 \subseteq [A:_R D] + Soc(R)$. Since $D$ is a Z-regular, then by Lemma 2.9, $(Soc(D) = Soc(R)D)$, and since $D$ is multiplication, then $A = [A:_R D]D$, and by Lemma 2.10 $rad_D(A) = \sqrt{[A:_R D]}D$, thus either $r_1T \subseteq rad_D(A) + Soc(D)$ or $r_2T \subseteq rad_D(A) + Soc(D)$ or $r_1r_2D \subseteq A + Soc(D)$. Therefore, by Proposition 2.11. $A$ is a weakly pseudo primary 2-absorbing sub-module of $D$. □

The following corollary is obtained as an immediate outcome of Proposition 2.12.

**Corollary 2.13**. Let $D$ be a cyclic Z-regular $R$-module and $A \subsetneq D$ be a sub-module of $D$. Then $A$ is a weakly pseudo primary 2-absorbing sub-module of $D$ if and only if $[A:_R D]$ is a weakly pseudo primary 2-absorbing ideal of $R$.

**Lemma 2.14**. [12] Let $D$ be a projective $R$-module, then $Soc(D) = Soc(R)D$.

**Proposition 2.15**. [5] Let $A \subsetneq D$ be a sub-module of $R$-module $D$, $A$ be a weakly pseudo primary 2-absorbing sub-module of an $R$-module $D$ if and only if $(0) \neq rQT \subseteq A$, for some $r \in R$ and for some ideal $Q$ of $R$ and sub-module $T$ of $D$, implies that either $rT \subseteq rad_D(A) + Soc(D)$ or $QT \subseteq rad_D(A) + Soc(D)$ or $rQ \subseteq [A + Soc(D):_R D]$.

**Proposition 2.16**. Let $D$ be a multiplication projective $R$-module and $A \subsetneq D$ be a sub-module of $D$. Then $A$ is a weakly pseudo primary 2-absorbing sub-module of $D$ if and only if $[A:_R D]$ is a weakly pseudo primary 2-absorbing ideal of $R$.

**Proof**. Assume that $(0) \neq dT_1T_2 \subseteq [A:_R D]$ for some sub-module $T_1, T_2$ of $R$, and $d \in D$ it follows that $(0) \neq dT_1T_2A \subseteq D$. But $D$ is a multiplication, then there exist two ideals $Q_1, Q_2$ in $R$ such that $T_1 = Q_1D$, $T_2 = Q_2D$. Implies that $(0) \neq dQ_1Q_2D \subseteq A$. Since $A$ is a weakly pseudo primary 2-absorbing sub-module of $D$. Then by Corollary 2.3, either $dQ_1D \subseteq rad_D(A) + Soc(D)$ or $dQ_2D \subseteq rad_D(A) + Soc(D)$ or $Q_1Q_2D \subseteq A + Soc(D)$. Since $D$ is a projective, then by Lemma 2.14, $(Soc(D) = Soc(R)D)$, and since $D$ is multiplication, then $A = [A:_R D]D$, and by Lemma 2.10 $rad_D(A) = \sqrt{[A:_R D]}D$, Thus either $dQ_1D \subseteq \sqrt{[A:_R D]}D + Soc(R)D$ or $dQ_2D \subseteq \sqrt{[A:_R D]}D + Soc(R)D$ or $Q_1Q_2D \subseteq [A:_R D]D + Soc(R)D$. That is either $dQ_1 \subseteq \sqrt{[A:_R D]} + Soc(R)$ or $dQ_2 \subseteq \sqrt{[A:_R D]} + Soc(R)$ or $Q_1Q_2 \subseteq [A:_R D] + Soc(R) = [[A:_R D] + Sco(R):_R R]$. Therefore, by Corollary 2.3. $[A:_R D]$ is a weakly pseudo primary 2-absorbing ideal of $R$.

Conversely, Assume that $(0) \neq rNT \subseteq A$ for some sub-module $T$ of $R$, and an ideal $N$ of $R$ and $r \in R$. But $D$ is a multiplication, then there exists an ideal $Q$ in $R$ such that $T = QD$. Implies that $(0) \neq rNQD \subseteq A$. that is $(0) \neq rNQ \in [A:_R D]$. Since $[A:_R D]$ is a weakly pseudo primary 2-absorbing ideal of $R$, then by Proposition 2.15, either $rN \subseteq \sqrt{[A:_R D]} + Soc(R)$ or $rQ \subseteq \sqrt{[A:_R D]} + Soc(R)$ or $NQ \subseteq [A:_R D] + Soc(R)$. Since $D$ is a projective, then by Lemma 2.14, $(Soc(D) = Soc(R)D)$, and since $D$ is multiplication, then $A = [A:_R D]D$, and by Lemma 2.10 $rad_D(A) = \sqrt{[A:_R D]}D$, thus either $rT \subseteq rad_D(A) + Soc(D)$ or $rT \subseteq rad_D(A) + Soc(D)$ or $NQ \subseteq A + Soc(D)$. Therefore, by Proposition 2.15. $A$ is a weakly pseudo primary 2-absorbing sub-module of $D$. □

The following corollary is obtained as an immediate outcome of Proposition 2.16.

**Corollary 2.17**. Let D be a cyclic and projective $R$-module and $A \subsetneq D$ be a sub-module of $D$. Then $A$ is a weakly pseudo primary 2-absorbing sub-module of $D$ if and only if $[A:_R D]$ is a weakly pseudo primary 2-absorbing ideal of $R$.

**Lemma 2.18**. [13] Let $D$ be a non-singular $R$-module, then $Soc(D) = Soc(R)D$.

**Proposition 2.19**. Let $D$ be a non-singular multiplication R-module and $A \subsetneq D$ be a sub-module of $D$. Then $A$ is a weakly pseudo primary 2-absorbing sub-module of $D$ if and only if $[A:_R D]$ is a weakly pseudo primary 2-absorbing ideal of $R$.

**Proof**. Assume that $(0) \neq T_1 T_2 d \subseteq [A:_R D]$ for some sub-modules $T_1, T_2$ of $R$ and $d \in D$, it follows that $(0) \neq T_1 T_2 dA \subseteq D$. But $D$ is a multiplication, then there exist three ideals $Q_1, Q_2$ in $R$ such that $T_1 = Q_1 D$ and $T_2 = Q_2 D$. Implies that $(0) \neq Q_1 Q_2 dD \subseteq A$. Since $A$ is a weakly pseudo primary 2-absorbing sub-module of $D$. Then by Corollary 2.3, either $Q_1 dD \subseteq rad_D(A) + Soc(D)$ or $Q_2 dD \subseteq rad_D(A) + Soc(D)$ or $Q_1 Q_2 D \subseteq A + Soc(D)$. Since $D$ is a non-singular, then by Lemma 2.18, $(Soc(D) = Soc(R)D)$, and since $D$ is multiplication, then $A = [A:_R D]D$, and by Lemma 2.10 $rad_D(A) = \sqrt{[A:_R D]}D$, thus either $Q_1 dD \subseteq \sqrt{[A:_R D]}D + Soc(R)D$ or $Q_2 dD \subseteq \sqrt{[A:_R D]}D + Soc(R)D$ or $Q_1 Q_2 D \subseteq [A:_R D]D + Soc(R)D$. That is either $Q_1 d \subseteq \sqrt{[A:_R D]} + Soc(R)$ or $Q_2 d \subseteq \sqrt{[A:_R D]} + Soc(R)$ or $Q_1 Q_2 \subseteq [A:_R D] + Soc(R) = [[A:_R D] + Sco(R):_R R]$. Therefore, by Corollary 2.3. $[A:_R D]$ is a weakly pseudo primary 2-absorbing ideal of $R$.

Conversely, Assume that $(0) \neq r_1 r_2 T \subseteq A$ for some sub-module $T$ of $R$, and $r_1, r_2 \in R$. But $D$ is a multiplication, then there exists an ideal $Q$ in $R$ such that $T = QD$. Implies that $(0) \neq r_1 r_2 QD \subseteq A$. That is $(0) \neq r_1 r_2 Q \in [A:_R D]$. Since $[A:_R D]$ is a weakly pseudo primary 2-absorbing ideal of $R$, Then we have by Proposition 2.1, either $r_1 Q \subseteq \sqrt{[A:_R D]} + Soc(R)$ or $r_2 Q \subseteq \sqrt{[A:_R D]} + Soc(R)$ or $r_1 r_2 \subseteq [A:_R D] + Soc(R)$. Since $D$ is a non-singular, then by Lemma 2.18, $(Soc(D) = Soc(R)D)$, and since $D$ is multiplication, then $A = [A:_R D]D$, and by Lemma 2.10 $rad_D(A) = \sqrt{[A:_R D]}D$, thus either $r_1 T \subseteq rad_D(A) + Soc(D)$ or $r_2 T \subseteq rad_D(A) + Soc(D)$ or $r_1 r_2 D \subseteq A + Soc(D)$. Therefore, by Proposition 2.1. $A$ is a weakly pseudo primary 2-absorbing sub-module of $D$. □

The following corollary is obtained as an immediate outcome of Proposition 2.19.

**Corollary 2.20.** Let $D$ be a cyclic and non-singular $R$-module and $A \subsetneq D$ be a sub-module of $D$. Then $A$ is a weakly pseudo primary 2-absorbing sub-module of $D$ if and only if $[A:_R D]$ is a weakly pseudo primary 2-absorbing ideal of $R$.

**Lemma 2.21.** [6] Let $D$ be a faithful multiplication $R$-module, then $Soc(D) = Soc(R)D$.

**Proposition 2.22.** Let $D$ be a multiplication faithful $R$-module, and $A \subsetneq D$ be a sub-module of $D$. Then $A$ is a weakly pseudo primary 2-absorbing sub-module of $D$ if and only if $[A:_R D]$ is a weakly pseudo primary 2-absorbing ideal of $R$.

**Proof.** The prove is similar to the prove of Proposition 2.19, and by Lemma 2.21. □

The following corollary is obtained as an immediate outcome of Proposition 2.22.

**Corollary 2.23.** Let D be a cyclic and faithful $R$-module and $A \subsetneq D$ be a sub-module of $D$. Then $A$ is a weakly pseudo primary 2-absorbing sub-module of $D$ if and only if $[A:_R D]$ is a weakly pseudo primary 2-absorbing ideal of $R$.

**Lemma 2.24.** [14] Let $D$ be a finitely generated multiplication $R$-module, and $Q_1, Q_2$ are ideals of $R$ then $Q_1 D \subseteq Q_2 D$, if and only if $Q_1 \subseteq Q_2 + ann_R(D)$.

**Lemma 2.25.** [15] Let $D$ be faithful multiplication $R$-module, then $rad_D(QD) = \sqrt{Q}D$ for an ideal $Q$ of $R$.

**Proposition 2.26.** Let D be a finitely generated Z-regular multiplication $R$-module and Q be an ideal of $R$ with $ann_R(D) \subseteq Q$. Then $Q$ is a weakly pseudo primary 2-absorbing ideal of $R$ if and only if $QD$ is a weakly pseudo primary 2-absorbing sub-module of $D$.

**Proof.** Assume that $(0) \neq T_1 T_2 T_3 \subseteq QD$ where $T_1, T_2, T_3$ are sub-module of $D$. Since $D$ is multiplication $R$-module, then $T_1 = Q_1 D$, $T_2 = Q_2 D$ and $T_3 = Q_3 D$ for some ideal $Q_1, Q_2, Q_3$ of $R$. that is $(0) \neq T_1 T_2 T_3 = Q_1 Q_2 Q_3 D \subseteq QD$. Since $D$ is a finitely generated multiplication $R$-module, then by Lemma 2.24 we have $(0) \neq Q_1 Q_2 Q_3 \subseteq Q + ann_R D$. Also, as $ann_R(D) \subseteq Q$, then $Q + ann_R D = Q$, that is $(0) \neq Q_1 Q_2 Q_3 \subseteq Q$. And as $Q$ is a weakly pseudo primary 2-absorbing ideal of $R$, then by Proposition 2.11 we obtain $Q_1 Q_3 \subseteq \sqrt{Q} + Soc(R)$ or $Q_2 Q_3 \subseteq \sqrt{Q} + Soc(R)$ or $Q_1 Q_2 \subseteq [Q + Soc(R):_R R] = Q + Soc(R)$. It follows either that $Q_1 Q_3 D \subseteq \sqrt{Q}D + Soc(R)D$ or $Q_2 Q_3 D \subseteq \sqrt{Q}D + Soc(R)D$ or $Q_1 Q_2 D \subseteq QD + Soc(R)D$. Since $D$ is a Z-regular then by Lemma 2.9, $(Soc(D) = Soc(R)D)$, and as $D$ is finitely multiplication by Lemma 2.25, $rad_D(QD) = \sqrt{Q}D$ it follows either that $T_1 T_3 \subseteq rad_D(QD) + Soc(D)$ or $T_1 T_3 \subseteq rad_D(QD) + Soc(D)$ or $T_1 T_2 \subseteq QD + Soc(D)$. Then by Proposition 2.11 we have $QD$ is a weakly pseudo primary 2-absorbing sub-module of $D$.

Conversely, let $QD$ be a weakly pseudo primary 2-absorbing sub-module of $D$, and let $(0) \neq Q_1 Q_2 Q_3 \subseteq Q$ for some ideals $Q_1, Q_2, Q_3$ in $R$ this implies, $(0) \neq Q_1 Q_2 Q_3 D \subseteq QD$, since $QD$ is a weakly pseudo primary 2-absorbing sub-module of $D$. Then by Proposition 2.11 either $Q_1 Q_3 D \subseteq rad_D(QD) + Soc(D)$ or $Q_2 Q_3 D \subseteq rad_D(QD) + Soc(D)$ or $Q_1 Q_2 D \subseteq QD + Soc(D)$ Now, as $D$ is Z-regular then by Lemma 2.9, $(Soc(D) = Soc(R)D)$, and since $D$ is finitely multiplication by Lemma 2.25, $rad_D(QD) = \sqrt{Q}D$. Then either $Q_1 Q_3 D \subseteq \sqrt{Q}D + Soc(R)D$ or $Q_2 Q_3 D \subseteq \sqrt{Q}D + Soc(R)D$ or $Q_1 Q_2 D \subseteq QD + Soc(R)D$. That is either $Q_1 Q_3 \subseteq \sqrt{Q} + Soc(R)$ or $Q_2 Q_3 \subseteq \sqrt{Q} +$ $Soc(R)$ or $Q_1Q_2 \subseteq QD + Soc(R)$. Then by Proposition 2.11. we have $Q$ is a weakly pseudo primary 2-absorbing ideal of $R$. □

The following corollary is obtained as an immediate outcome of Proposition 2.26.

**Corollary 2.27**. Let $D$ be a cyclic finitely generated $Z$-regular $R$-module and $Q$ be an ideal of $R$ with $ann_R(D) \subseteq Q$. Then $Q$ is a weakly pseudo primary 2-absorbing ideal of $R$ if and only if $QD$ is a weakly pseudo primary 2-absorbing sub-module of $D$.

**Proposition 2.28**. Let $D$ be a finitely generated multiplication projective $R$-module and $Q$ be an ideal of $R$ with $ann_R(D) \subseteq Q$. Then $Q$ is a weakly pseudo primary 2-absorbing ideal of R if and only if $QD$ is a weakly pseudo primary 2-absorbing sub-module of $D$.

**Proof**. Assume that $(0) \neq T_1T_2T_3 \subseteq QD$ where $T_1, T_2, T_3$ are sub-modules of $D$. Since $D$ is multiplication $R$-module, then $T_1 = Q_1D$, $T_2 = Q_2D$ and $T_3 = Q_3D$ for some ideals $Q_1, Q_2, Q_3$ of $R$, so $(0) \neq T_1T_2T_3 = Q_1Q_2Q_3D \subseteq QD$. Since $D$ is a finitely generated multiplication $R$-module, then by Lemma 2.24 we have $(0) \neq Q_1Q_2Q_3 \subseteq Q + ann_RD$. Since $ann_R(D) \subseteq Q$, then $Q + ann_RD = Q$, that is $(0) \neq Q_1Q_2Q_3 \subseteq Q$. And as $Q$ is a weakly pseudo primary 2-absorbing ideal of $R$, then by Proposition 2.11 $Q_1Q_3 \subseteq \sqrt{Q} + Soc(R)$ or $Q_2Q_3 \subseteq \sqrt{Q} + Soc(R)$ or $Q_1Q_2 \subseteq [Q + Soc(R):_R R] = Q + Soc(R)$. It follows either that $Q_1Q_3D \subseteq \sqrt{Q}D + Soc(R)D$ or $Q_2Q_3D \subseteq \sqrt{Q}D + Soc(R)D$ or $Q_1Q_2D \subseteq QD + Soc(R)D$. Since $D$ is projective then by Lemma 2.14, $(Soc(D) = Soc(R)D)$, and since $D$ is finitely multiplication by Lemma 2.25, $rad_D(QD) = \sqrt{Q}D$. It follows either that $T_1T_3 \subseteq rad_D(QD) + Soc(D)$ or $T_1T_3 \subseteq rad_D(QD) + Soc(D)$ or $T_1T_2 \subseteq QD + Soc(D)$. Then by Proposition 2.11 we have $QD$ is a weakly pseudo primary 2-absorbing sub-module of $D$.

Conversely, let $QD$ be a weakly pseudo primary 2-absorbing sub-module of $D$, and let $(0) \neq Q_1Q_2Q_3 \subseteq Q$ for some ideals $Q_1.Q_2.Q_3$ in $R$. Which implies that, $(0) \neq Q_1Q_2Q_3D \subseteq QD$, as $QD$ is a weakly pseudo primary 2-absorbing sub-module of $D$. Then by Proposition 2.11 either $Q_1Q_3D \subseteq rad_D(QD) + Soc(D)$ or $Q_2Q_3D \subseteq rad_D(QD) + Soc(D)$ or $Q_1Q_2D \subseteq QD + Soc(D)$ Since $D$ is projective then by Lemma 2.14, $(Soc(D) = Soc(R)D)$, and as $D$ is finitely multiplication by Lemma 2.25, $rad_D(QD) = \sqrt{Q}D$. Then either $Q_1Q_3D \subseteq \sqrt{Q}D + Soc(R)D)$ or $Q_2Q_3D \subseteq \sqrt{Q}D + Soc(R)D$ or $Q_1Q_2D \subseteq QD + Soc(R)D$. That is either $Q_1Q_3 \subseteq \sqrt{Q} + Soc(R)$ or $Q_2Q_3 \subseteq \sqrt{Q} + Soc(R)$ or $Q_1Q_2 \subseteq QD + Soc(R)$. Then by Proposition 2.11. we have $Q$ is a weakly pseudo primary 2-absorbing ideal of $R$. □

The following corollary is obtained as an immediate outcome of Proposition 2.28.

**Corollary 2.29**. Let $D$ be a cyclic finitely generated projective $R$-module and $Q$ be an ideal of $R$ with $ann_R(D) \subseteq Q$. Then $Q$ is a weakly pseudo primary 2-absorbing ideal of $R$ if and only if $QD$ is a weakly pseudo primary 2-absorbing sub-module of $D$.

**Proposition 2.30**. Let $D$ be a finitely generated multiplication non-singular $R$-module and $Q$ be an ideal of $R$ with $ann_R(D) \subseteq Q$. Then $Q$ is a weakly pseudo primary 2-absorbing ideal of $R$ if and only if QD is a weakly pseudo primary 2-absorbing sub-module of $D$.

**Proof.** Assume that $(0) \neq T_1 T_2 T_3 \subseteq AD$ where $T_1, T_2, T_3$ are sub-modules of $D$. Since $D$ is multiplication $R$-module, then $T_1 = Q_1 D$, $T_2 = Q_2 D$ and $T_3 = Q_3 D$ for some ideals $Q_1, Q_2, Q_3$ of $R$. That is $(0) \neq T_1 T_2 T_3 = Q_1 Q_2 Q_3 D \subseteq QD$. Since $D$ is a finitely generated multiplication $R$-module, then by Lemma 2.24 we have $(0) \neq Q_1 Q_2 Q_3 \subseteq Q + ann_R D$. Since $ann_R(D) \subseteq Q$, then $Q + ann_R D = Q$, that is $(0) \neq Q_1 Q_2 Q_3 \subseteq Q$. Since $A$ is a weakly pseudo primary 2-absorbing ideal of $R$, then by Proposition 2.11 $Q_1 Q_3 \subseteq \sqrt{Q} + Soc(R)$ or $Q_2 Q_3 \subseteq \sqrt{Q} + Soc(R)$ or $Q_1 Q_2 \subseteq [Q + Soc(R) :_R R] = Q + Soc(R)$. It follows either that $Q_1 Q_3 D \subseteq \sqrt{Q}D + Soc(R)D$ or $Q_2 Q_3 D \subseteq \sqrt{Q}D + Soc(R)D$ or $Q_1 Q_2 D \subseteq QD + Soc(R)D$. Since $D$ is non-singular then by Lemma 2.18, $(Soc(D) = Soc(R)D)$, and since $D$ is finitely multiplication by Lemma 2.25, $rad_D(QD) = \sqrt{Q}D$. It follows either that $T_1 T_3 \subseteq rad_D(QD) + Soc(D)$ or $T_1 T_3 \subseteq rad_D(QD) + Soc(D)$ or $T_1 T_2 \subseteq QD + Soc(D)$. Then by Proposition 2.11 we have $QD$ is a weakly pseudo primary 2-absorbing sub-module of $D$.

Conversely, let $QD$ be a weakly pseudo primary 2-absorbing sub-module of $D$, and let $(0) \neq Q_1 Q_2 Q_3 \subseteq Q$ for some ideals $Q_1. Q_2. Q_3$ in $R$. Implies that, $(0) \neq Q_1 Q_2 Q_3 D \subseteq QD$, since $QD$ is a weakly pseudo primary 2-absorbing sub-module of $D$. Then by Proposition 2.11 either $Q_1 Q_3 D \subseteq rad_D(QD) + Soc(D)$ or $Q_2 Q_3 D \subseteq rad_D(QD) + Soc(D)$ or $Q_1 Q_2 D \subseteq QD + Soc(D)$ Since $D$ is non-singular then by Lemma 2.18, $(Soc(D) = Soc(R)D)$, and since $D$ is finitely multiplication by Lemma 2.25, $rad_D(QD) = \sqrt{Q}D$. Then either $Q_1 Q_3 D \subseteq \sqrt{Q}D + Soc(R)D)$ or $Q_2 Q_3 D \subseteq \sqrt{Q}D + Soc(R)D$ or $Q_1 Q_2 D \subseteq QD + Soc(R)D$. That is either $Q_1 Q_3 \subseteq \sqrt{Q} + Soc(R)$ or $Q_2 Q_3 \subseteq \sqrt{Q} + Soc(R)$ or $Q_1 Q_2 \subseteq QD + Soc(R)$. Then by Proposition 2.11, we have $Q$ is a weakly pseudo primary 2-absorbing ideal of $R$. □

The following corollary is obtained as an immediate outcome of Proposition 2.30.

**Corollary 2.31.** Let $D$ be a cyclic finitely generated non-singular $R$-module and $Q$ be an ideal of $R$ with $ann_R(D) \subseteq Q$. Then Q is a weakly pseudo primary 2-absorbing ideal of $R$ if and only if QD is a weakly pseudo primary 2-absorbing sub-module of $D$.

**Proposition 2.32.** Let $D$ be a faithful finitely generated multiplication $R$-module and $Q$ be an ideal of $R$ with $ann_R(D) \subseteq Q$. Then $Q$ is a weakly pseudo primary 2-absorbing ideal of $R$ if and only if QD is a weakly pseudo primary 2-absorbing sub-module of $D$.

**Proof.** The prove is similar to the prove of Proposition 2.30, and by Lemma 2.21. □

The following corollary is obtained as an immediate outcome of Proposition 32.

**Corollary 2.33.** Let $D$ be a cyclic faithful finitely generated $R$-module and $Q$ be an ideal of $R$ with $ann_R(D) \subseteq Q$. Then $Q$ is a weakly pseudo primary 2-absorbing ideal of $R$ if and only if $QD$ is a weakly pseudo primary 2-absorbing sub-module of $D$.

**Lemma 2.34.** [16] Let $D$ be a multiplication $R$-module, then $D$ is cancellation if and only if $D$ is a faithful finitely generated.

**Proposition 2.35.** Let D be a faithful finitely generated multiplication $R$-module, and $A \subsetneq D$ be a sub-module of $D$. The below statements are considered to be equal:

  (a)  $A$ is a weakly pseudo primary 2-absorbing sub-module of $D$.

(b) $[A:_R D]$ is a weakly pseudo primary 2-absorbing ideal of $R$.

(c) $A = QD$ for some $Q$ be a weakly pseudo primary 2-absorbing ideal of $R$.

**Proof.** $(a) \Leftrightarrow (b)$ Follows by Proposition 2. 22.

$(b) \Rightarrow (c)$ Assume that $[A:_R D]$ is a weakly pseudo primary 2-absorbing ideal of $R$. since $D$ is multiplication then $A = [A:_R D]D$. Put $Q = [A:_R D]$ is a weakly pseudo primary 2-absorbing ideal of $R$ and $A = QD$.

$(c) \Rightarrow (b)$ Assume that $A = QD$ for some weakly pseudo primary 2-absorbing ideal of $R$. Since $D$ is multiplication, then $A = [A:_R D]D = QD$. Since $D$ is faithful finitely generated $R$-module, then by Lemma 2.34, $D$ is a cancellation, hence $[A:_R D] = Q$ is a weakly pseudo primary 2-absorbing ideal of $R$. □

The following result is an immediate outcome of a Proposition 2.35.

**Corollary 2.36.** Let $D$ be a faithful cyclic $R$-module and $A \subsetneq D$ be a sub-module of $D$. The below statements are considered to be equal:

(a) $A$ is a weakly pseudo primary 2-absorbing sub-module of $D$.

(b) $[A:_R D]$ is a weakly pseudo primary 2-absorbing ideal of $R$.

(c) $A = QD$ for some $Q$ be a weakly pseudo primary 2-absorbing ideal of $R$.

**Lemma 2.37.** [16] Let $D$ be a multiplication $R$-module, then $D$ be finitely generated if and only if $D$ is a weak cancellation.

**Proposition 2.38.** Let $D$ be Z-regular and finitely generated multiplication $R$-module and $A \subsetneq D$ be a sub-module of $D$. With $ann_R(D) \subseteq [A:_R D]$. The below statements are considered to be equal:

(a) $A$ is a weakly pseudo primary 2-absorbing sub-module of $D$.

(b) $[A:_R D]$ is a weakly pseudo primary 2-absorbing ideal of $R$.

(c) $A = QD$ for some $Q$ be a weakly pseudo primary 2-absorbing ideal of $R$ with $ann_R(D) \subseteq Q$.

**Proof.** $(a) \Leftrightarrow (b)$ Follows by Proposition 2.12

$(b) \Rightarrow (c)$ Let $D$ be sub-module of $D$, then $A = [A:_R D]D$ (for $D$ is multiplication). put $Q = [A:_R D]D$ implies that $Q$ is a weakly pseudo primary 2-absorbing ideal of $R$ with $ann_R(D) = [0:D] \subseteq [A:_R D] = Q$, that is $ann_R(D) \subseteq Q$

$(c) \Rightarrow (b)$ Assume that $A = QD$ for some $Q$ be a weakly pseudo primary 2-absorbing ideal of $R$, with $ann_R(D) \subseteq Q$. Since $D$ is multiplication then $A = [A:_R D]D = QD$. Since $D$ is faithful finitely generated $R$-module, then by Lemma 2.37. $D$ is a weak cancellation, therefore $[A:_R D] + ann_R(D) = Q + ann_R(D)$. But $ann_R(D) \subseteq Q$ and $ann_R(D) \subseteq [A:_R D]$, that is $[A:_R D] = Q$. since $Q$ is a weakly pseudo primary 2-absorbing ideal of $R$. Then $[A:_R D]$ is a weakly pseudo primary 2-absorbing ideal of $R$. □

The following result is an immediate outcome of Proposition 2.38.

**Corollary 2.39.** Let $D$ be cyclic $Z$-regular $R$-module and $A \subsetneq D$ be a sub-module of $D$. With $ann_R(D) \subseteq [A:_R D]$. The below statements are considered to be equal:

(a) $A$ is a weakly pseudo primary 2-absorbing sub-module of $D$.

(b) $[A:_R D]$ is a weakly pseudo primary 2-absorbing ideal of $R$.

(c) $A = QD$ for some $Q$ be a weakly pseudo primary 2-absorbing ideal of $R$ with $ann_R(D) \subseteq Q$.

**Proposition 2.40.** Let $D$ be projective and finitely generated multiplication projective $R$-module and $A \subsetneq D$ be a sub-module of $D$. With $ann_R(D) \subseteq [A:_R D]$. The below statements are considered to be equal:

(a) $A$ is a weakly pseudo primary 2-absorbing sub-module of $D$.

(b) $[A:_R D]$ is a weakly pseudo primary 2-absorbing ideal of $R$.

(c) $A = QD$ for some Q be a weakly pseudo primary 2-absorbing ideal of $R$ with $ann_R(D) \subseteq Q$.

**Proof.** $(a) \Leftrightarrow (b)$ Follows by Proposition 2.12

$(b) \Leftrightarrow (c)$ Follows by Proposition 2.38. □

The following result is an immediate outcome of Proposition 2.40.

**Corollary 2.41.** Let $D$ be a cyclic projective $R$-module and $A \subsetneq D$ be a sub-module of $D$. With $ann_R(D) \subseteq [A:_R D]$ The below statements are considered to be equal:

(a) $A$ is a weakly pseudo primary 2-absorbing sub-module of $D$.

(b) $[A:_R D]$ is a weakly pseudo primary 2-absorbing ideal of $R$.

(c) $A = QD$ for some $Q$ be a weakly pseudo primary 2-absorbing ideal of $R$ with $ann_R(D) \subseteq Q$.

**Proposition 2.42.** Let $D$ be non-singular and finitely generated multiplication $R$-module and $A \subsetneq D$ be a sub-module of $D$. With $ann_R(D) \subseteq [A:_R D]$ The below statements are considered to be equal:

(a) $A$ is a weakly pseudo primary 2-absorbing sub-module of $D$.

(b) $[A:_R D]$ is a weakly pseudo primary 2-absorbing ideal of $R$.

(c) $A = QD$ for some $Q$ be a weakly pseudo primary 2-absorbing ideal of $R$ with $ann_R(D) \subseteq Q$.

**Proof.** $(a) \Leftrightarrow (b)$ Follows by Proposition 2.19.

$(b) \Leftrightarrow (c)$ Follows by Proposition 2.38. □

The following result is obtained as an immediate outcome of a Proposition 2.42.

**Corollary 2.43.** Let $D$ be a cyclic non-singular $R$-module and $A \subsetneq D$ be a sub-module of $D$. With $ann_R(D) \subseteq [A:_R D]$ The below statements are considered to be equal:

(a) $A$ is a weakly pseudo primary 2-absorbing sub-module of $D$.

(b) $[A:_R D]$ is a weakly pseudo primary 2-absorbing ideal of $R$.

(c) $A = QD$ for some $Q$ be a weakly pseudo primary 2-absorbing ideal of $R$ with $ann_R(D) \subseteq Q$.

## Bibliography

[1] Y. A. Darani and F. Soheilnia, "2-absorbing and weakly 2-absorbing submodules," *Thai J. Math,* vol. 9, no. 3, pp. 577-584, 2011.

[2] A. Badawi, "On 2-absorbing ideals of commutative rings," *Bulletin of the Australian Mathematical,* vol. 75, no. 3, pp. 417-429, 2007.

[3] H. Mostafanasab, E. Yetkin, U. Tekir and A. Darani, "On 2-absorbing primary submodules of modules over commutative rings," *Analele ştiinţifice ale Universităţii" Ovidius" Constanţa. Seria Matematică,* vol. 24, no. 1, pp. 335-351, 2016.

[4] O. Abdulla and H. Mohammadali, "Pseudo Primary-2-Absorbing sub-modules and Some Related Concepts," *Ibn Al-Haitham Journal for Pure and Applied Sciences,* vol. 32, no. 3, 2019.

[5] O. H. Taha and M. A. Salih, "Weakly Pseudo Primary 2-Absorbing Submodules," *Al-Bahir Journal for Engineering and Pure Sciences,* vol. 5, no. 1, 2024.

[6] A. A. Tuganbaev, "Multiplication modules," *Journal of Mathematical Sciences,* vol. 123, no. 2, pp. 3839-3905, 2004.

[7] J. Zelmanowitz, "Semiprime modules with maximum conditions," *Journal of Algebra,* vol. 25, no. 3, pp. 554-574, 1973.

[8]  C. Faith, Algebra: rings, modules and categories I, Springer Science & Business Media, 2012.

[9]  A. Goldie, Torsion-free modules and rings, vol. 1, University of Chicago Press, 1972, p. 268–287.

[10] T-Y. Lam, Exercises in modules and rings, Springer Science & Business Media, 2019.

[11] R. Wisbauer and R. Wisbauer, Foundations of module and ring theory, Routledge, 2018.

[12] R. McCasland and M. Moore, "On radicals of submodules," *Communications in Algebra,* vol. 19, no. 5, pp. 1327-1341, 1991.

[13] L. Rowen, Ring theory: 83, Academic Press, 2012.

[14] C. W. CHOI and E. S. KIM, "Some remarks on multiplication modules," *Tamkang Journal of Mathematics,* vol. 24, no. 3, pp. 309-316, 1993.

[15] M. M. Ali, "Invertibility of multiplication modules II," *New Zealand J. Math,* vol. 39, pp. 45-64, 2009.

[16] S. Ali, "On cancellation modules," Diss. M. Sc. Thesis, University of Baghdad, 1993.